\documentclass[11pt,reqno]{amsart}
\usepackage{amsmath}
\usepackage{amsthm}
\usepackage{amssymb}
\usepackage{enumerate}
\usepackage{amscd}
\usepackage{color}
\usepackage{pb-diagram}
\usepackage{graphicx}

\theoremstyle{plain}

\newtheorem*{prop*}{Proposition}
\theoremstyle{definition}

\theoremstyle{remark}

\numberwithin{equation}{section}

\newcommand{\vs}{\vspace{0.2cm}}
\newcommand{\vvs}{\vspace{0.4cm}}

\begin{document}
\title[Real K\"ahler Submanifolds in Codimension Four]{ An extension theorem for Real K\"ahler Submanifolds in Codimension Four}
\author{Jinwen Yan}
\address{Center for Mathematical Sciences\\
                    Zhejiang University\\
                    Hangzhou, 310027 China}
\email{yimkingman@gmail.com}
\author{Fangyang zheng}

\address{Center for Mathematical Sciences\\Zhejiang University \\
Hangzhou, 310027 China, and  Department of Mathematics\\
                    The Ohio State University\\
                    231 West 18th Avenue, Columbus, OH 43210\\}
\email{zheng@math.ohio-state.edu}

\begin{abstract} In this article, we prove a K\"ahler extension theorem for real K\"ahler submanifolds of codimension $4$ and rank at least $5$. Our main theorem  states that such a manifold is a holomorphic hypersurface in another real K\"ahler submanifold of codimension $2$. This generalizes a result of Dajczer and Gromoll  in 1997 which states that any real K\"ahler submanifolds of codimension $3$ and rank at least $4$ admits a K\"ahler extension.
\end{abstract}
\maketitle

\section{Introduction}

\vvs

Submanifold theory, and especially the study of Riemannian submanifolds in Euclidean spaces, have been a classic subarea in differential geometry. The Nash embedding theorem \cite{Nash} guarantees that any complete Riemannian manifold can be isometrically embedded into an Euclidean space. There are lots of important development in submanifold theory. At the risk of omitting many, we will just mention two recent such examples. One is the work of Hongwei Xu and his collaborators (see \cite{Xu}, \cite{XG} and \cite{XZ}) generalizing the famous Differentiable Sphere Theorem of Brendle and Schoen (\cite{BS1}, \cite{BS2}) to the submanifold case, thus obtaining the optimal pinching constant. The other one is the very recent work by F. Marques and A. Neves \cite{MN}, solving the long-standing Willmore conjecture.

\vs

However, in the special case when the submanifold happens to be K\"ahler, the research is relatively few and sporadic, and the state of knowledge is still rather primitive in our opinion. We will call a K\"ahler manifold isometrically embedded in a real Euclidean space a {\em real K\"ahler Euclidean submanifold,} or {\em real K\"ahler submanifold} for short.  That is, we have an isometric embedding $f: M^n \rightarrow {\mathbb R}^{2n+p}$ from a  K\"ahler manifold $M^n$ of complex dimension $n$ into the real Euclidean space.

\vs

Ideally, since $M^n$ is equipped with a complex structure, one would like the embedding $f$ to be both isometric and holomorphic. However, the thesis of Calabi \cite{Calabi} in 1950's showed us that very few K\"ahler metrics can be isometrically and holomorphically embedded in a complex Euclidean space or other complex space forms. He actually precisely characterized all such metrics. So to study generic K\"ahler manifolds in the extrinsic setting, one has to abandon the holomorphicity assumption on the embedding, and only assume it to be isometric.

\vs

For a real K\"ahler submanifold $f: M^n \rightarrow {\mathbb R}^{2n+p}$, the K\"ahlerness of $M^n$ imposes  strong restrictions and made it very sensitive to its codimension. For instance, when $p=1$, namely, when $M^n$ is a hypersurface, the result of Florit and the second named author in \cite{FZ-hypersurface} states that, when $M^n$ is also assumed to be complete, $f$ must be the product of $g$ with the identity map of ${\mathbb C}^{n-1}$, where $g:\Sigma \rightarrow {\mathbb R}^3$ is the isometric embedding of a complete surface, which is always K\"ahler. In other words, surfaces in ${\mathbb R}^3$ are essentially the only real K\"ahler submanifolds in codimension one. In contrast, there are all kinds of real hypersurfaces in Euclidean spaces.

\vs

In codimension two, the situation is also well-studied and fully understood. In the minimal case, it was analyzed in details by Dajczer and Gromoll (see \cite{DG2}, \cite{DR} and the references therein), and in the non-minimal case, it was classified by Florit and the second named author \cite{FZ-codim2}. In codimension three, the work of Dajczer and Gromoll \cite{DG97} showed that, unless the submanifold $M^n$ is a holomorphic hypersurface  of a real K\"ahler submanifold of codimension $1$, its rank has to be less than or equal to $3$, the codimension of $M^n$.

\vs

Recall that the {\em rank} of a real K\"ahler submanifold $f: M^n \rightarrow {\mathbb R}^{2n+p}$ at $x\in M$ is defined to be $n-\nu_0$, with $\nu_0$ the complex dimension of $\Delta_0=\Delta \cap J\Delta$, which is the $J$-invariant part of the kernel $\Delta $ of the second fundamental form of $f$. Of course these spaces may not have constant dimensions on $M$. But if we let $U$ be the open subset where $\Delta_0$ takes the minimum (thus constant) dimension, then $r$ will be constant in $U$. Outside the closure of $U$, $M$ will be a real K\"ahler submanifold with smaller rank. In general, by restricting to an open dense subset $U'$ of $M$, we can always assume that in each connected component $U$ of $U'$,   $\Delta$ and $\Delta_0$ take constant dimensions and form distributions. Note that the leaves of $\Delta $ ($\Delta_0$) are totally geodesic (complex) submanifolds in $M^n$. They are actually open subset of (parallel translation of) linear subspaces in the ambient Euclidean space. We might need to further reduce $U'$ later, but the conclusions we will draw will always be valid in each connected component of an open dense subset of $M$.

\vs

The main purpose of this paper is to show that the result of Dajczer and Gromoll in \cite{DG97} can be extended to the  codimension $4$ case. To be precise, we will prove the following

\vvs

\noindent {\bf Main Theorem.} {\em Let $f: M^n \rightarrow {\mathbb R}^{2n+4}$ be a real K\"ahler submanifold with rank $r>4$ everywhere. Then there exists an open dense subset $U'\subset M$ such that for each connected component $U$ of $U'$, the restriction $f|_U$ has a K\"ahler extension, namely, there exists a real K\"ahler submanifold  $h: Q^{n+1} \rightarrow {\mathbb R}^{2n+4}$ of codimension $2$, and a holomorphic embedding $\sigma : U \rightarrow Q^{n+1}$, such that $f|_U=h \circ \sigma $. Furthermore, when $f$ is minimal, one can choose $h$ to be minimal as well.}

\vvs

Note that if $h$ is minimal, $f$ has to be minimal. In general, the extension $h$ might not be unique. But as we shall see from the proof, there is always a `canonical' extension, unless $f$ itself is a holomorphic isometric embedding into ${\mathbb C}^{n+2}$.

\vs

This result can be regarded as an extension of a phenomenon discovered by Dajczer \cite{D} and Dajczer-Gromoll \cite{DG97}, in codimension two and three, respectively. In \cite{D}, Dajczer proved that, for any codimension two real K\"ahler submanifold, if its rank is greater than $2$, then in any connected component $U$ of an open dense subset of $M$, the restriction $f|_U$ is a holomorphic embedding into ${\mathbb R}^{2n+2}\cong {\mathbb C}^{n+1}$. This is an important discovery. In codimension three, Dajczer and Gromoll proved in 1997 (see \cite{DG97}) that if real K\"ahler submanifold of dimension three has rank greater than $3$, then there exists an open dense subset $U'\subseteq M$ such that in each connected component $U$ of $U'$, $f|_U$ has a K\"ahler extension into a real K\"ahler submanifold $Q^{n+1}$ of codimension one.

\vs

Note that for the results in \cite{D} and \cite{DG97}, assumptions were made on the {\em relative nullity} $\nu$, namely, the (real) dimension of the kernel $\Delta $ of the second fundamental form $\alpha_f$. Since $\Delta_0\subseteq \Delta$, we have $2\nu_0\leq \nu$ hence $\nu \geq 2n-2r$, with $r$ the rank. In \cite{D}, the assumption was $\nu <2n-4$, which implies $r>2$. In \cite{DG97}, the assumption was $\nu <2n-6$, which implies $r>3$. Even though their assumptions were slightly stronger, it is easy to see that their arguments can be extended to the cases when assumptions are made on the ranks.

\vs

We suspect that similar phenomenon will persist in higher codimensions as well, namely, the rank $r$ should be controlled by the codimension $p$ in a certain way, unless the manifold is a complex submanifold of another real K\"ahler submanifold of a smaller codimension. We will explore the higher codimensional cases elsewhere, but in here we will just state a conjecture which says that, for $p\leq 11$, the words ``controlled by" in the above sentence should mean that the rank is no greater than the codimension, namely, $r\leq p$. In other words,

\vvs

\noindent {\bf Conjecture.} {\em Let $f: M^n \rightarrow {\mathbb R}^{2n+p}$ be a real K\"ahler submanifold with rank $r>p$ everywhere. If $p\leq 11$, then there exists an open dense subset $U'\subset M$ such that for each connected component $U$ of $U'$, the restriction $f|_U$ has a K\"ahler extension, namely, there exists a real K\"ahler submanifold  $h: Q^{n+s} \rightarrow {\mathbb R}^{2n+p}$ of codimension $p-2s<p$, and a holomorphic embedding $\sigma : U \rightarrow Q^{n+s}$, such that $f|_U=h \circ \sigma $. }

\vvs

Note that the main theorem, together with results of \cite{D} and \cite{DG97}, confirms the conjecture for $p\leq 4$. (When $p=1$, one always has $r\leq 1$).

\vvs

\noindent {\bf Acknowledgement:} We would like to take this opportunity to thank a few people who helped us in our study. First, we are very grateful to Marcos Dajczer for his inspiring papers on the subject of real K\"ahler submanifolds, which opened the way to the investigation of this under-explored territory in submanifold theory. The second named author would like to thank his former collaborators Luis Florit and Wing San Hui. The present work is a continuation of these earlier joint works. Finally, we would also like to thank CMS of Zhejiang University which provides an ideal research environment for mathematicians, and in particular to Hongwei Xu for his warm hospitality and numerous stimulating conversations.

\vvs

\vvs

\section{Preliminaries}

\vvs

In this section, we shall collect some known results in the literature that will be needed in the proof of our theorem. We will also fix some notations and terminologies that will be used later.

\vs

In this paper, unless specified otherwise, we will always assume that $M$ is a real K\"ahler submanifold of complex dimension $n$ and codimension $p$, with $f$ the isometric embedding from $M$ into ${\mathbb R}^{2n+p}$. At any $x\in M$, let $\Delta $ be the kernel of the second fundamental form $\alpha_f$ of $f$, and $ \Delta_0 =\Delta \cap J\Delta $ the $J$-invariant part of $\Delta$. The rank $r$ is defined to be $n-\nu_0$, where $2\nu_0$ is the real dimension of $\Delta_0$. We always have $\nu \geq 2n-2r$, where $\nu =\mbox{dim}(\Delta )$ is the relative nullity.

\vs

The results in this paper are  local in nature, and we will from time to time reduce from $M$ into an open dense subset of it, to make various subspaces in the tangent or normal bundle taking constant dimensions and forming subbundles.

\vs

For $x\in M$, we will denote by $T\cong {\mathbb R}^{2n}$ the real tangent space $T_xM$, $N=T_xM^{\perp }\cong {\mathbb R}^p$ the normal space, and by $V\cong {\mathbb C}^n$ the space of all type $(1,0)$ complex tangent vectors at $x$, namely, $V\oplus \overline{V} \cong T\otimes_{\mathbb R}{\mathbb C}$. Extend the second fundamental form $\alpha_f: T\times T \rightarrow N$  linearly over ${\mathbb C}$, we will denote its $(1,1)$ and $(2,0)$ components by $H$ and $S$, respective:
$H: V\otimes \overline{V} \rightarrow N_{\mathbb C}$, and $S:  V\otimes V \rightarrow N_{\mathbb C}$,
where $N_{\mathbb C}=N\otimes_{\mathbb R}{\mathbb C}$.

\vs

As observed in  \cite{FHZ}, the K\"ahlerness of $M$ implies that the Hermitian bilinear form $H$ and the symmetric bilinear form $S$ satisfy the following symmetry conditions:
\begin{eqnarray}
\langle H_{X\overline{Y}}, H_{Z\overline{W}}\rangle & = & \langle H_{Z\overline{Y}}, H_{X\overline{W}}\rangle  \\
\langle H_{X\overline{Y}}, S_{ZW}\rangle & = & \langle H_{Z\overline{Y}}, S_{XW}\rangle  \\
\langle S_{XY}, S_{ZW}\rangle & = & \langle S_{ZY}, S_{XW}\rangle
\end{eqnarray}
for any $X,Y,Z,W \in V$.

\vs

We notice that $H$ and $S$ together carry all the information of $\alpha_f$. Also, by $(2.1)$, we get
$$ | \sum_{i=1}^n  H_{i\overline{i}} |^2 = \sum_{i,j=1}^n |H_{i\overline{j}}|^2 $$
for any unitary frame $\{ e_1, \ldots , e_n\}$ of $V$. Here we wrote $H_{i\overline{j}}$ for $H_{e_i\overline{e_j}}$. So $H\equiv 0$ if and only if the trace of $H$, which is (a multiple of) the mean curvature of $f$, vanishes. So $f$ is minimal when and only when $H=0$.

\vs

Note that for $\Delta = \mbox{ker}(\alpha_f)$, its $J$-invariant part $\Delta_0=\Delta \cap J\Delta$ corresponds to a complex subspace $D\subseteq V$ with complex dimension $\nu_0$, and $D$ is exactly the intersection of the kernels of $H$ and $S$. Let $V'$ be the orthogonal complement of $D$ in $V$. We have $V=D\oplus V'$ and $V'\cong {\mathbb C}^r$, where $r=n-\nu_0$ is the rank of $M^n$. $D$ (or $\Delta$) is contained in the kernel of the curvature tensor of $M$, and the leaves of the foliation $D$ are totally geodesic, flat complex submanifolds in $M$. They are actually open subset of ${\mathbb C}^{n-r}$, embedded linearly (i.e., as parallel translation of linear subspace) in ${\mathbb R}^{2n+p}$. So in a way, the rank $r$ of $M$  is like the essential (complex) dimension of $M$, even though in general $M$ might not be isometric to the product space (i.e., the leaves of $D$ might not be parallel to each other).

\vs

For any $\eta \in N$, the shape operator $A_{\eta}$ is defined by $\langle A_{\eta }u,v\rangle =\langle \alpha_f(u,v), \eta \rangle$ for any $u, v\in T$. It is self-adjoint. For convenience, we will also denote by $A^{\eta }$ the {\em shape form}, which is defined by $A^{\eta}_{uv} = \langle A_{\eta }(u), v\rangle = \langle \alpha_f(u,v), \eta \rangle $. It is the component  of the second fundamental form in the $\eta$-direction.

\vs

Let $\{ e_1, \ldots , e_n\}$ be a basis of $V$. For each $1\leq i\leq n$, write
$$ e_i=\frac{1}{\sqrt{2}} (\varepsilon _i -\sqrt{-1} \varepsilon_{n+i}).$$
Then under the basis $\{ \varepsilon_1, \ldots , \varepsilon_{2n}\}$ of $T$, $A^{\eta}$ will take the form
\begin{eqnarray}
 A^{\eta} & = & \left( \begin{array}{ll} \mbox{Re}(H^{\eta })+\mbox{Re}(S^{\eta }) & \mbox{Im}(H^{\eta })- \mbox{Im}(S^{\eta }) \\  - \mbox{Im}(H^{\eta })-\mbox{Im}(S^{\eta }) & \mbox{Re}(H^{\eta })-\mbox{Re}(S^{\eta })\end{array} \right)
\end{eqnarray}
where $H^{\eta }=\langle H_{i\overline{j}}, \eta \rangle$ and $S^{\eta }=\langle S_{ij}, \eta \rangle$. Note that under any tangent frame $\{\varepsilon_1, \ldots , \varepsilon_{2n}\}$, the shape operator $A_{\eta }$ and the shape form $A^{\eta }$ are related by
$$ A_{\eta} (\varepsilon_i) = \sum_{j=1}^{2n} (A^{\eta }g^{-1})_{ij}\varepsilon_j =\sum_{j,k=1}^{2n} A^{\eta }_{ik}g^{kj}\varepsilon_j
$$
where $A^{\eta}_{ij}=A^{\eta}_{\varepsilon_i\varepsilon_j}$, $g_{ij}=\langle \varepsilon_i, \varepsilon_j\rangle$, and $(g^{ij})$ is the inverse matrix of $(g_{ij})$.

\vvs

Next let us recall the Codazzi equation:
\begin{eqnarray}
\nabla_u(A_{\xi }v)-\nabla_v(A_{\xi }u)-A_{ \nabla^{\perp }_u\xi }v+ A_{ \nabla^{\perp }_v\xi }u -A_{\xi }[u,v]=0
\end{eqnarray}
for any vector fields $u$, $v$ on $M$ and normal section $\xi$. For any type $(1,0)$ tangent vector $X$ and any (possibly complexified) normal vector $\xi$, let us denote by
\begin{eqnarray}
 A_{\xi }X=H_{\xi }X + S_{\xi } X
\end{eqnarray}
the decomposition of $A_{\xi }X $  into its $(1,0)$ part and $(0,1)$ part. This give us operators $H_{\xi }$ and $S_{\xi }$ which are determined by
$$ H_{\xi }X=\sum_{i=1}^n  H^{\xi }_{X\overline{i}} e_i, \ \ \ \ S_{\xi }X=\sum_{i=1}^n S^{\xi }_{Xi} \overline{e_i} . $$
under any unitary frame $\{ e_1, \ldots , e_n\}$ of $V$. Note that $H_{\xi }(V)\subseteq V$ and $H_{\xi }(\overline{V})\subseteq \overline{V}$; \ while $S_{\xi }(V)\subseteq \overline{V}$ and $S_{\xi }(\overline{V})\subseteq V$. Extend the Codazzi equation linearly to all complexified tangent vectors, and by taking the $(1,0)$ and $(0,1)$ parts in (2.5), we get
\begin{eqnarray}
 \nabla_X(H_{\xi }Y)-\nabla_Y(H_{\xi }X)-H_{ \nabla^{\perp }_X\xi }Y + H_{ \nabla^{\perp }_Y\xi }X -H_{\xi }[X,Y] &= &0 \ \ \ \ \\
 \nabla_{X}(S_{\xi }Y)-\nabla_{Y}(S_{\xi }X)-S_{ \nabla^{\perp }_X\xi }Y+ S_{ \nabla^{\perp }_Y\xi }X -S_{\xi }[X,Y] &=&0
\end{eqnarray}
and
$$
\nabla_{\overline{Y}}(S_{\xi }X)-S_{\nabla^{\perp }_{\overline{Y}}\xi }X-S_{\xi }(\nabla_{\overline{Y}}X)=
 \nabla_{X} (H_{\xi }{\overline{Y}})-H_{\nabla^{\perp }_X\xi }\overline{Y}-H_{\xi }(\nabla_X\overline{Y})
$$
for any type $(1,0)$ vector fields $X$, $Y$ on $M$ and any normal field $\xi$. In particular, in the minimal case, namely, when $H=0$, we have
\begin{eqnarray}
S_{\nabla^{\perp }_{\overline{Y}}\xi }X = \nabla_{\overline{Y}}(S_{\xi }X)-S_{\xi }(\nabla_{\overline{Y}}X),  \ \ \ \ \ \mbox{when} \ H=0
\end{eqnarray}
for any $\xi$ in $N$ and any $X$, $Y$ in $V$.

\vvs

\vvs

\section{The Algebraic Lemma}

\vvs

In this paper, we shall be primarily interested in the case when $p=4$ and $r>4$, although some of the arguments work in general cases as well. Our first objective is to show that at a generic point $x$ in $M^n$, the second fundamental form takes a rather special form. First, let us introduce the following

\vvs

\noindent {\bf Definition.} {\em Let $V\cong {\mathbb C}^n$ and $N\cong {\mathbb R}^p$ be equipped with inner products, and let $H$, $S$ be respectively Hermitian or symmetric bilinear map from $V$ into $N_{\mathbb C}=N\otimes {\mathbb C}$ satisfying the symmetry conditions (2.1)-(2.3). Let $E$ be a subspace of $N$. An {\bf almost complex structure} $J$ on $E$ is an isometry from $E$ onto itself, such that $J^2=-I$, and for any $\eta \in E$, $H^{\eta}=0$ and $S^{J\eta }=-\sqrt{-1}S^{\eta }$ holds. }

\vvs

Here we wrote $H^{\eta} =\langle H, \eta\rangle$ and $S^{\eta} =\langle S, \eta\rangle$. Note that $E$ is necessarily even dimensional, and the condition on $J$ is equivalent to $A_{J\eta }= JA_{\eta}$ for any $\eta \in E$, where $A_{\eta}$ is the shape operator, related to the shape form $A^{\eta}$ by the metric on $T\cong V$, which in turn is related to $H^{\eta}$ and $S^{\eta }$ by (2.4).

\vs

We will assume that the dimension $p$ of $N$ is the smallest, namely, for any $\eta \neq 0$ in $N$, either $H^{\eta}$ or $S^{\eta }$ is not zero. This is equivalent to $A_{\eta }\neq 0$ for any $\eta \neq 0$ in $N$. Note that under this assumption, the almost complex structure on any subspace $E$ of $N$, if exists, must be unique. To see this, suppose $J$ and $J'$ are both almost complex structures on $E\subseteq N$. Then for any $\eta \in E$, we have $H^{\eta }=0$ and $S^{J\eta }=-\sqrt{-1}S^{\eta} = S^{J'\eta }$, so $S^{J\eta -J'\eta }=0$. So if $J\neq J'$, then by (2.4) there will be $\eta \neq 0$ in $E$ such that $A_{\eta }=0$, contradicting our assumption that $p$ is the smallest.

\vs

As a consequence of this uniqueness, we know that if $E_1$, $E_2$ are both subspaces of $N$ admitting almost complex structures, then both $E_1\cap E_2$ and $E_1+E_2$ also admit almost complex structure. So there is always a (unique) maximal subspace $E$ in $N$, possibly trivial, that is equipped with an almost complex structure. We will call this subspace $E$ the {\em complex part of} $N$.

\vs

Let $E'$ be the orthogonal complement of the complex part $E$ in $N$, and write $S'=\langle S, E'\rangle $. Then by the definition of the almost complex structure, we know $S'$ again satisfies (2.3). Also, if $S^{\eta }$ has rank at most $1$, then in $\{\eta\}^{\perp}$, $S$ also satisfies (2.3). Our main goal in this section is to prove the following

\vvs

\noindent {\bf Algebraic Lemma.} {\em Let $V\cong {\mathbb C}^r$, $N\cong {\mathbb R}^4$ be equipped with inner products, and let $H$, $S$  be respectively Hermitian or symmetric bilinear forms from $V$ into $N_{\mathbb C}$ satisfying symmetry conditions (2,1)-(2.3). We assume that $\mbox{ker}(H)\cap \mbox{ker}(S)=0$ and $r>4$. Then $N$ has non-trivial complex part. That is, either $N$ itself or a $2$-dimensional subspace $E$ in it admits an almost complex structure. Furthermore, in the latter case we have
$$\mbox{dim}(\mbox{ker}(H) \cap \mbox{ker}(S'))\geq r-2,$$
 where $S'=\langle S, E'\rangle $ and $E'$ is the orthogonal complement of $E$ in $N$. }

\vvs

\noindent {\em Proof:} Since $H$ is Hermitian, its image space is in the form $N'_{\mathbb C}=N'\otimes {\mathbb C}$ for some real linear subspace $N'\subseteq N$. Let $N=N'\oplus N''$ be the orthogonal decomposition and write $H=(H',H'')$ and $S=(S',S'')$ under this decomposition. We have $H''=0$ by definition. Denote by $p'$, $q=4-p'$ the dimension of $N'$, $N''$, respectively.

\vs

Let $V_0$ be the kernel of $H$, and $V=V_0\oplus V_1$ the orthogonal decomposition. Write $r_i=\mbox{dim}_{\mathbb C}V_i$ for $i=0,1$. Note that for any $X\in V_0$, $H_{X\overline{\ast }} =0$, so by (2.2), we know that $\langle S_{XY},H_{\ast \overline{\ast }}\rangle =0$ thus $S'_{XY}=0$, for any $Y\in V$. Hence $V_0\subseteq \mbox{ker}(S')$.

\vs

From the discussion in  \cite{FHZ}, we know that $r_1\leq p'$, and the equality case would imply that $H'$ and $S'$ can be simultaneously diagonalized. In particular, $p'=4$ cannot happen, since $r\geq 5$. Similarly, $p'=3$ cannot happen, either. This is  because in this case the rank of $S'$ is at most $r_1\leq 3$. The fact $r\geq 5$ and the symmetry condition (2.3) would make $S''$, thus $S$, having a zero-eigenvector within $V_0$, contradicting the fact that $\mbox{ker}(H)\cap \mbox{ker}(S)=0$ in $V$. So we have $p'\leq 2$.

\vs

If $p'=2$, then $r_1$ is necessarily $2$, and we are in the diagonal situation. That is, we will have orthonormal basis $\{\xi_1, \xi_2\}$ of $N'$ and basis $\{ e_1, e_2\}$ of $V_1$ such that $V_0=\mbox{ker}(H)\cap \mbox{ker}(S')$, and along $V_1$, the matrices $H^1$, $H^2$, $S^1$, and $S^2$ are respectively
\begin{eqnarray*}
 \left( \begin{array}{ll} 1 & 0 \\ 0 & 0  \end{array} \right)  ; \ \   \left( \begin{array}{ll} 0 & 0 \\ 0 & 1 \end{array} \right) ;\ \ \left( \begin{array}{ll} \ast & 0 \\ 0 & 0  \end{array} \right)  ; \ \ \left( \begin{array}{ll} 0 & 0 \\ 0 & \ast  \end{array} \right)
\end{eqnarray*}

Notice that both $S^1$ and $S^2$ have rank $\leq 1$, so the symmetric bilinear form $S''$ from $V$ into $N''\cong {\mathbb R}^2$ satisfies (2.3) as well. Its kernel cannot overlap with $V_0$, so its rank is at least $3$. By Lemma 1 below, we know that $N''$ admits an almost complex structure.

\vs

If $p'=1$, then $r_1=1$ necessarily, so $V_1$ is one-dimensional and both $H'$ and $S'$ are zero in the codimension one subspace $V_0$ of $V$. Since $S'$ is a matrix of rank $\leq 1$, the remaining part $S''$ will satisfy (2.3) and its rank is at least $4$. So by Lemma 1 below, $N''$ contains a $2$-dimensional subspace $E$ which admits an almost complex structure. Let $0\neq \eta \in N''$ be perpendicular to $E$. Then $S^{\eta }$ again satisfies (2.3), so its rank is at most $1$. Putting $\eta $ together with $N'$ to form the space $E'$, we know that the common kernel of $H$ and $S$ on $E'$ has dimension at least $r-2$.

\vs

Finally, when $p'=0$, we are left with $S$ from $V$ into $N={\mathbb R}^4$ satisfying (2.3) and with rank at least $5$. So by Lemma 1 below, we know that either $N$ itself admits an almost complex structure, or it contains a $2$-dimensional subspace $E$ which does. Let $E'=E^{\perp }$ in $N$. Since $S'=\langle S,E'\rangle$ also satisfies (2.3), if it does not admit an almost complex structure, then by Lemma 1 is must have rank less than or equal to $2$, namely, $\mbox{dim(ker}(S))\geq r-2$. This completes the proof of the Algebraic Lemma. \qed

\vvs

\noindent {\bf Lemma 1.} {\em Let $V\cong {\mathbb C}^r$ and $N\cong {\mathbb R}^p$ be equipped with inner products, write $N_{\mathbb C}=N\otimes {\mathbb C}$. Let $S: V\times V \rightarrow N_{\mathbb C}$ be a symmetric bilinear map, satisfying (2.3) and with  $\mbox{ker}(S)=0$. If $p\leq 4$ and $r>p$, then there exists $X,Y\in V$ such that $S_{XY}\neq 0$ and $\langle S_{XY}, S_{ZW}\rangle =0$ for any $Z,W\in V$. In other words, $N$ always has nontrivial complex part. }

\vvs

\noindent {\em Proof:} The $p=2$ case is due to Dajczer in \cite{D}, and the $p=3$ case is due to Dajczer and Gromoll \cite{DG97}, even though their notations are quite different from here.  We will just prove the $p=4$ case here, since the same argument would work for the $p=2$ and $p=3$ cases as well. Without loss of generality, we may assume that $r=5$ (as when $r>5$, we can just apply the result to any $5$-dimensional subspace of $V$).

\vs

For $X\in V$, consider the linear map $\phi_X:V\rightarrow N_{\mathbb C}$ sending $Y$ to $S_{XY}$. Denote by $K_X$ the kernel of $\phi_X$, and $k_X$ its complex dimension. Since $V\cong {\mathbb C}^5$, $N_{\mathbb C} \cong {\mathbb C}^4$, and $\mbox{ker}(S)=0$, we have $1\leq k_X\leq 4$.

\vs

Let $k$ be the minimum of $k_X$ for all $X\in V$, and denote by $V_0$ be the open dense subset of $V$ consisting of all $X$ with $k_X=k$. We will also write $m=5-k$. It is the dimension of the image of $\phi_X$ and is also between $1$ and $4$. Notice that the set $\Sigma =\{ X\in V\mid S_{XX}=0\}$ is the intersection of four quadratic hypersurfaces in $V$, so $V_0'=V_0\setminus \Sigma $ is still open dense in $V$.

\vs

Fix any $X\in V_0'$. Let $\{ e_1, \ldots , e_5\}$ be a basis of $V$ such that $e_1=X$, $\{ e_{m+1}, \ldots , e_5\}$ forms a basis of $K_X$. Again we will write $S_{ij}$ for $S_{e_ie_j}$. $\{ S_{11}, \ldots , S_{1m}\}$ forms a basis of the image space $P=\phi_X(V)$. We will denote by $Q$ the subspace of $N_{\mathbb C}$ spanned by $S_{i\alpha }$ for all $1\leq i\leq 5$ and all $m< \alpha \leq 5$. That is, $Q=S(K_X\times V)$. Since $S_{1\alpha }=0$, the symmetry condition (2.3) implies that $\langle P, Q\rangle =0$.

\vs

We claim that $Q\subseteq P$. Assume otherwise. Then there will some $m<\alpha \leq 5$ and some $1\leq i\leq 5$, such that $S_{i\alpha }$ is not contained in $P$. Consider the vector $Y=e_1+\lambda e_i$ for a sufficiently small $\lambda$. Then $S_{Y\alpha }=\lambda S_{i\alpha }$, and we have
\begin{eqnarray*}
 S_{Y1} \wedge \cdots \wedge S_{Ym} \wedge S_{Y\alpha }= \lambda (S_{11}\wedge \cdots \wedge S_{1m}\wedge S_{i\alpha } +O(\lambda))
\end{eqnarray*}
whose leading term is not zero. So for a sufficiently small value of $\lambda$, the image of $\phi_Y$ has dimension bigger than $m$, a contradiction. This proves that $Q\subseteq P$. Note that $Q\neq 0$ since $\mbox{ker}(S)=0$.

\vs

When $m=1$, $Q=P$, so $0\neq S_{11}\in P=Q$ satisfies $\langle S_{11}, S_{ij}\rangle =0$ for any $i$, $j$. If $m=2$, then since we can take $e_2\in V_0'$ also, both $K_1$ and $K_2$ are of codimension $2$, thus there will be $0\neq Z\in K_1\cap K_2$. Take $W$ such that $S_{ZW}\neq 0$, then $S_{ZW}\in Q$, and $\langle S_{ZW}, S_{22}\rangle =0$, hence $\langle S_{ZW}, S_{ij}\rangle =0$ for any $i$, $j$. On the other hand, since $\langle P,Q\rangle =0$, $P$ is contained in the orthogonal complement of $\overline{Q}$ in $N_{\mathbb C}$, so $m\leq 3$. From now on, we will assume that $m=3$.

\vs

Note that if there are $\alpha , \beta \in \{ 4,5\}$ such that $S_{\alpha \beta }\neq 0$, then since $\langle Q, Q\rangle =0$, by (2.3), we would have
$$ \langle S_{\alpha \beta }, S_{ij}\rangle = \langle S_{\alpha i }, S_{\beta j}\rangle =0 $$
for any $i,j\leq 3$. So $S_{\alpha \beta }$ will give us the proof of the lemma. In other words, if for some $X\in V_0'$ we have $S(K_X\times K_X)\neq 0$, then any non-zero element $S_{ZW}$ in this subspace would satisfy $\langle S_{ZW}, S_{ij}\rangle =0$ for all $i$, $j$. So we may further assume that $S(K_X\times K_X)= 0$ for all $X\in V_0'$. We claim that this will not be possible at all, thus completing the proof of the lemma.

\vs

Since $V_0'$ is open dense in $V$. We may assume that $e_2$, $e_3$ are in $V_0'$ also. Consider their kernels $K_2$ and $K_3$. If they are both equal to $K_1$, then $e_4$ will be in the kernel of $S$, a contradiction. So we must have one of them, say $K_2$, not equal to $K_1$. Since $Q$ has dimension $1$, $S_{24}$ and $S_{25}$ are proportional to each other. Replace $\{ e_4, e_5\}$  by another basis of $K_1$ if necessary, we may assume that $S_{24}=0$. On the other hand, since $K_2\neq K_1$, we may replace $e_3$ by another vector in $K_2$. So $K_2=\mbox{span}\{ e_3, e_4\}$. Since $e_2\in V_0'$, we know that $S(K_2\times K_2)=0$ (unless the lemma holds). However, this means $S_{34}=S_{44}=0$. But we already have $S_{14}=S_{54}=0$ since $e_4\in K_1$, hence $e_4\in \mbox{ker}(S)$, a contradiction once again. This finishes the proof of the lemma. \qed

\vvs

\vvs

\section{The Extension Theorems}

\vvs

Now let us consider a real K\"ahler submanifold $f: M^n \rightarrow {\mathbb R}^{2n+4}$ of codimension $4$. Reduce $M$ to a connected component $U$ of an open dense subset $U'$ of $M$ if necessary, we may assume that both $\Delta $ and $\Delta_0$ are of constant dimensions and are distributions. We will also assume that at any $x\in M$, the shape operator $A_{\xi }\neq 0$ for any $\xi \neq 0$. Note that the vanishing of some shape operator everywhere would mean that the codimension can be reduced. By the algebraic lemma proved in the previous section, we know that either the entire normal bundle $N$ or a rank two subbundle $E\subseteq N$ admits an almost complex structure.

\vs

We will call an almost complex structure $J$ on $E$ {\em admissible} if
\begin{eqnarray}
J (\nabla^{\perp }_v\xi)^E = (\nabla^{\perp }_vJ\xi)^E
\end{eqnarray}
holds for any $\xi \in E$ and any vector field $v$ in $M$. Here $(W)^E $ stands for the $E$ component of $W$.

\vs

Notice that in the case when $E$ has rank $2$, any almost complex structure $J$ on $E$ is automatically admissible: let $\{\xi_1,\xi_2\}$ be a local orthonormal frame of $E$ with $\xi_2=J\xi_1$. Equation (4.1) reduces to
$$ J( \langle \nabla^{\perp }\xi_1, \xi_2\rangle \xi_2 )= \langle \nabla^{\perp }\xi_2 , \xi_1\rangle \xi_1  ,$$
or equivalently
$$ \langle \nabla^{\perp }\xi_1, \xi_2\rangle = - \langle \nabla^{\perp }\xi_2 , \xi_1\rangle ,$$
which always holds.

\vs

In the case when $N$ itself admits an admissible almost complex structure $J$, our goal is to show that $M^n$ is actually a holomorphic submanifold in ${\mathbb C}^{n+2}$. We have the following:

\vvs

\noindent {\bf Theorem 1. } {\em Let $f: M^n \rightarrow {\mathbb R}^{2n+4}$  be a real K\"ahler submanifold whose normal bundle  admits an admissible almost complex structure. Then there exists an isometric identification $\sigma : {\mathbb R}^{2n+4}\cong {\mathbb C}^{n+2}$ such that $\sigma \circ f$ is a holomorphic isometric embedding.}

\vvs

We will prove this theorem at the end of this section.

\vs

In the case of a rank two subbundle $E$ of $N$ admitting  an almost complex structure, we would like to show that $M^n$ is a complex submanifold of another complex manifold $Q^{n+1}$, and $Q^{n+1}$ is a codimension two real K\"ahler submanifold of which $M$ is the restriction. We will call such a $Q^{n+1}$ a {\em K\"ahler extension} of $M^n$. To prove this extension theorem, we will need to know more information about the behavior of the second fundamental form beyond the existence of the almost complex structure on $E$. It turns out that what is needed here is the following data:

\vvs

\noindent {\bf Definition.} {\em A {\em developable ruling in $E\oplus T$}  is a rank two subbundle $L$ of $E\oplus T$,  such that $L+T=E\oplus T$ and  $\langle \widetilde{\nabla }L, E'\rangle =0$ along $M$. Here $T$ is the tangent bundle of $M$, $E'$ is the orthogonal complement of $E$ in the normal bundle $N$, and $\widetilde{\nabla }$ is the covariant differentiation of the ambient Euclidean metric.}

\vvs

Note that the subbundle $L$ is necessarily transversal to $T$, but in general not contained in $N$. We will prove the following extension theorem:

\vvs

\noindent {\bf Theorem 2. } {\em Let $f: M^n \rightarrow {\mathbb R}^{2n+4}$ be a real K\"ahler submanifold. If there is a rank two subbundle $E$ of the normal bundle $N$, an almost complex structure $J$ on $E$, and a developable ruling $L$ in $E\oplus T$. Then there exists a real K\"ahler submanifold $h: Q^{n+1}\rightarrow {\mathbb R}^{2n+4}$ and a holomorphic embedding $\sigma : M^n \rightarrow Q^{n+1}$ such that $f=h\circ \sigma $. }

\vvs

\noindent {\em Proof:} Let $z=(z_1, \ldots , z_n)$ be a local holomorphic coordinate in $M$ and $\{ \xi_1, \ldots , \xi_4\}$ be an orthonormal frame of $N$, such that $\{ \xi_1,\xi_2\}$ spans $E'$ and $\{ \xi_3,\xi_4\}$ spans $E$. Write $P=E\oplus T$. Since $L+T=P$, there will be a local frame of $L$ given by
$$ \eta_1=\xi_3-v_1, \ \ \ \ \eta_2=\xi_4-v_2 $$
where $v_1$ and $v_2$ are real vector fields of $M$. Since $\langle \widetilde{\nabla }L, E'\rangle =0$, we know that
\begin{eqnarray}
\widetilde{\nabla }_v\eta_i \ \in \ P = L+T
\end{eqnarray}
for $i=1$, $2$ and for any vector field $v$ in $M$.

\vs

Let $B\subseteq {\mathbb C}$ be a sufficiently small disc and $t=t_1+\sqrt{-1}t_2$ be the coordinate. Define a $(2n+2)$-dimensional submanifold $h: Q\rightarrow {\mathbb R}^{2n+4}$ by
$$ h(z,t)=f(z)+t_1\eta_1(z)+t_2\eta_2(z) $$
Since $L$ is transversal to $T$, for sufficiently small values of $|t|$ the map $h$ is an embedding. $Q$ is ruled along the directions of $L$. By $(4.2)$, the bundle $E'$, which is the normal bundle of $Q$, is constant along each leave of $L$, thus $Q$ is a developable submanifold (meaning that its tangent space is constant along each ruling). Along the submanifold $M$ of $Q$, the restriction of the tangent bundle $TQ|_M$ is simply $P=L+T$. Since $P=E\oplus T$, and we have almost complex structure $J$ on both $T$ and $E$, we can take their direct sum to get an almost complex structure on $P$. Now take parallel translation along leaves of $L$, we get an almost complex structure on $TQ$. We will denote this almost complex structure on $TQ$ again by $J$.

\vs

To show that $Q$ is a K\"ahler manifold under the restriction of the Euclidean metric, it suffices to show that $\widehat{\nabla}J=0$ on $Q$, where $\widehat{\nabla}$ is the connection on $Q$, namely, the $Q$-component  of $\widetilde{\nabla}$. That is, we just need to show that
\begin{eqnarray} \widehat{\nabla}_Z(JW) & = & J(\widehat{\nabla}_ZW)
\end{eqnarray}
holds for any two vector fields $Z$ and $W$ in $Q$. Since $TQ$ is the parallel translation in ${\mathbb R}^{2n+4}$ of $TQ|_M=P$ along the leaves of $L$, and $J$ is also defined by parallel translation along leaves of $L$, we just need to verify the above condition at points in $M$ and with $Z$ tangent to $M$. If $W$ is also tangent to $M$, then the above equation holds in the tangential component of $M$, since $M$ is K\"ahler. For the normal components, since we are only concerned within $Q$, it means that we just need to verify that for the $\xi_3$ and $\xi_4$ directions, namely:
\begin{eqnarray*} \langle \widehat{\nabla}_Z( J W) ,\xi_i \rangle & = & \langle J(\widehat{\nabla}_Z W) , \xi_i \rangle
\end{eqnarray*}
for $i=3$ and $4$ where $Z$ and $W$ are vector fields in $M$. It is equivalent to
\begin{eqnarray} J A_{\xi_i} & = & A_{J\xi_i}
\end{eqnarray}
for $i=3$ and $4$.
Since $H^{\xi_3}=H^{\xi_4}=0$, $S^{\xi_3}=\sqrt{-1}S^{\xi_4}$, so by (2.4) we get
$$ JA_{\xi_3}=\left( \begin{array}{cc} 0&-1\\1&0 \end{array} \right)   \left( \begin{array}{cc} R_3&-I_3\\-I_3&-R_3 \end{array} \right) =\left( \begin{array}{cc} I_3&R_3\\R_3&-I_3 \end{array} \right) = A_{\xi_4}.$$
Here we wrote $S^{\xi_3}=R_3+\sqrt{-1}I_3$ and $S^{\xi_4}=R_4+\sqrt{-1}I_4$, so $R_3=-I_4$ and $I_3=R_4$. Recall that we have defined $J$ on $E$ by $J\xi_3=\xi_4$ and $J\xi_4=-\xi_3$. So (4.4) holds.

\vs

 Now we are left with the case where $Z$ is a tangent vector field of $M$ and $W$ is a section of $E$, since $P=E\oplus T$. By the linearity of $J$ and the Leibniz formula, we just need to check this for $W=\xi_3$ and $W=\xi_4$. Namely,
\begin{eqnarray}
\widehat{\nabla}_Z(\xi_4)&=& J(\widehat{\nabla}_Z\xi_3)
\end{eqnarray}
for any tangent vector field $Z$ in $M$.  First let us compare the tangential components on both sides. It reduces once again to (4.4). For the normal components in (4.5), notice that $\widehat{\nabla}$ is just the $TQ$ component of $\widetilde{\nabla}$, so we have
\begin{eqnarray*}
 (\widehat{\nabla}_Z\xi_3)^{\perp }& = & \langle \widehat{\nabla}_Z\xi_3, \xi_4\rangle \xi_4 \ = \ \langle \nabla^{\perp }_Z\xi_3, \xi_4\rangle \xi_4 \ = \ -\langle \xi_3, \nabla^{\perp }_Z\xi_4\rangle \xi_4 \\
 (\widehat{\nabla}_Z\xi_4)^{\perp }& = &\langle \widehat{\nabla}_Z\xi_4, \xi_3\rangle \xi_3 \ =\ \langle \nabla^{\perp }_Z\xi_4, \xi_3\rangle \xi_3
 \end{eqnarray*}
 So $(J\widehat{\nabla}_Z\xi_3)^{\perp} = J((\widehat{\nabla}_Z\xi_3)^{\perp}) =(\widehat{\nabla}_Z\xi_4)^{\perp }$. This proves the K\"ahlerness of the codimension $2$ submanifold $Q$ in the Euclidean space. The holomorphicity of $M$ in $Q$ is obvious, since we defined our $J$ on $Q$ in such a way that its restriction on $M$ comes from the complex structure. This completes the proof of Theorem 2. \qed

\vvs

For the K\"ahler extension $h$ obtained in Theorem 2, clearly, if $h$ is minimal, then $f$ is necessarily minimal. Conversely, when $f$ is minimal, we would like to know when will $h$ be minimal. We have the following

\vvs

\noindent {\bf Theorem 3.} {\em Let $f$, $(E,J)$ and $L$ be as in Theorem 2, and let $h$ be the K\"ahler extension of $f$ obtained by $L$. Suppose $f$ is minimal, then $h$ is minimal if and only  $(v_2-Jv_1)\in \mbox{ker}(A_{\xi_1})\cap \mbox{ker}(A_{\xi_2})$. Here $\{ \xi_1 , \ldots , \xi_4\}$ is an orthonormal frame of $N$, with $\{ \xi_3, \xi_4\}$ a frame of $E$,  $\xi_4=J\xi_3$, and $v_1, v_2\in T$ are determined (uniquely) by the condition that $\{ \xi_3-v_1, \xi_4-v_2\}$ spans $L$.    }

\vvs

\noindent {\em Proof:} Note that $\xi_1$ and $\xi_2$ span the normal bundle of $Q$ in ${\mathbb R}^{2n+4}$, and $h$ is minimal if and only if its $H=0$, or equivalently, $J\hat{A}_{\xi_{\alpha }}=\hat{A}_{\xi_{\alpha }}J$ for $\alpha =1 $ and $2$, where $J$ is the almost complex structure of $Q$ and $\hat{A}$ is the shape operator of $Q$. That is, for $1\leq \alpha \leq 2$ and any vector fields $Z$, $W$ on $Q$,
\begin{eqnarray*}
 \langle J\hat{A}_{\xi_{\alpha }}Z, W\rangle =  \langle \hat{A}_{\xi_{\alpha }}J Z, W\rangle ,
 \end{eqnarray*}
 or equivalently,
\begin{eqnarray}
 - \langle \widetilde{\nabla}_Z J W , \xi_{\alpha } \rangle =  \langle \widetilde{\nabla}_{JZ} W , \xi_{\alpha } \rangle  .
 \end{eqnarray}
By the construction of $h$, $TQ$ is the parallel translate of of $TQ|_M$ along the leaves of $L$, and $J$ and both $\xi_{\alpha }$ are parallel along each leaf of $L$, so we just need to check (4.6) at points in $M$, and for $Z$ a vector field in $M$.

\vs

Since $TQ|_M=E\oplus T$, we just need to verify (4.6) for $W$ being a vector field in $M$ and a section of $E$. In the former case, (4.6) is just the minimality of $f$. While when $W$ is a section of $E$, (4.6) becomes
\begin{eqnarray}
  \langle  J W , \widetilde{\nabla}_Z\xi_{\alpha } \rangle =  - \langle  W , \widetilde{\nabla}_{JZ}\xi_{\alpha } \rangle
 \end{eqnarray}
for each $\alpha =1,2$. Clearly, we just need to verify (4.7) for $W=\xi_3$.

\vs

Now suppose that $\xi_3-v_1$ and $\xi_4-v_2$  span $L$, and $\xi_4=J\xi_3$. Note that since $L$ is transversal to $T$, the map $\pi|_L: L\rightarrow E$ is bijective. Here $\pi$ is the projection map from $E\oplus T$ onto $E$. So $v_1$, $v_2$ are uniquely determined by the choice of $\{ \xi_3,\xi_4\}$. By the definition of developable ruling, we know that $\langle \widetilde{\nabla}\xi_{\alpha }, L\rangle =0$, so
 \begin{eqnarray*}
 \langle  \xi_4 , \widetilde{\nabla}_Z\xi_{\alpha } \rangle &=& \langle  v_2 , \widetilde{\nabla}_Z\xi_{\alpha } \rangle \ = \ \langle A_{\xi_{\alpha }}(v_2), Z\rangle , \ \  \mbox{and} \\
\langle  \xi_3 , \widetilde{\nabla}_{JZ}\xi_{\alpha } \rangle &= & \langle  v_1 , \widetilde{\nabla}_{JZ}\xi_{\alpha } \rangle \ = \ \langle A_{\xi_{\alpha }} (v_1), JZ\rangle \ = \ \langle A_{\xi_{\alpha }} (J v_1), Z\rangle
\end{eqnarray*}
Note that in the last equality we used the minimality of $M$, namely, we always have $JA=-AJ$. Plug these two equalities into (4.7) for $W=\xi_3$, we get
$$ \langle A_{\xi_{\alpha }}(v_2-Jv_1) , Z\rangle =0 $$
for any vector field $Z$ in $M$, that is
\begin{eqnarray}
 A_{\xi_{\alpha }}(v_2-J v_1)=0, \ \ \ \ \alpha = 1, 2.
 \end{eqnarray}
So when $f$ is minimal, $h$ will be minimal if and only $v_2-Jv_1$ belongs to $\mbox{ker}(A_{\xi_1})\cap \mbox{ker}(A_{\xi_2})$, which is the real subspace of $T$ corresponding to $\mbox{ker}(S')$ in $V$. Here $S'=(S^1,S^2)$. This completes the proof of Theorem 3. \qed

\vvs

\noindent {\em Remark:} \ Let us denote by $\pi : E\oplus T \rightarrow E$ the projection map, and by $\tau : E \rightarrow L$ the inverse of the restriction map $\pi|_L: L\rightarrow E$. Then the condition stated in Theorem 3 can be rephrased as
$$ \tau (J\eta)-J\tau (\eta) \ \in \ \mbox{ker}(A_{\xi_1})\cap \mbox{ker}(A_{\xi_2}) $$
for any $\eta$ in $E$. Here $\{ \xi_1, \xi_2\}$ is a basis of $E'$, the orthogonal complement of $E$ in $N$.

\vs

Now let us prove Theorem 1 stated at the beginning of this section.

\vvs

\noindent {\em Proof of Theorem 1:} Note that in this case, the ambient Euclidean space is automatically a developable submanifold (of itself) over $M$, with fibers of the normal bundle $N$ as rulings leaves. Define an almost complex structure $J$ on $T\oplus N$ by taking the direct sum of the almost complex structure of $M$ with the given one on $N$, and use parallel translation along leaves of $N$ to push it to a small tubular neighborhood $\Omega $ of $M$, we get an almost complex structure $J$ on the open subset $\Omega $ of ${\mathbb R}^{2n+4}$. $J$ is clearly an isometry. One can see that $\widetilde{\nabla }J=0$ just like in the proof of Theorem 2, with the help of $(4.1)$. So this $J$ comes from an isometric identification ${\mathbb R}^{2n+4}\cong {\mathbb C}^{n+2}$ and $M$ becomes a complex submanifold with complex codimension $2$. This completes the proof of Theorem 1. \qed

\vvs

\vvs

\section{The Proof of the Main Theorem}

\vvs

In this section, we will prove the main theorem. For $x\in M$, let us denote by $N_0(x)$ the subspace of $N_x$ consisting of all $\eta$ with $A_{\eta }=0$. Note that the presence of normal directions in which the shape operator vanishes  would mean that the codimension can be reduced (see \cite{Spivak}, Prop. 24). In the interior part $U_0$ of the set where $N_0\neq 0$, there will be open dense subset of $U_0$, such that within each connected component of it the submanifold $M$ will be real K\"ahler  submanifold with smaller codimensions. Since the main theorem is known in codimension three or less, in the following, we will assume that

\vvs

{\em $N_0=0$ everywhere in $M$. That is, $A_{\eta }\neq 0$ for any $\eta \neq 0$.}

\vvs

First let us consider the non-minimal case, in other words, we restrict ourselves to the open subset of $M$ in which $H\neq 0$, if that set is non-empty. Since $r\geq 5$, we know that the image of $H$ is either $1$ or $2$ dimensional. In the open subset $U_2$ where $H$ has $2$-dimensional image space $E'$, there are exactly two directions, perpendicular to each other, in which $H$ has rank $1$. Let $\xi_1$ and $\xi_2$ be the unit vectors in those two directions, they are unique up to $\pm 1$ and interchange. In this case, as a consequence of (2.2), $S^{\xi_1}$ and $S^{\xi_2}$ can be diagonalized accordingly.

\vs

In the open subset $M\setminus \overline{U_2}$, the image of $H$ is $1$-dimensional, and we will let $\xi_1$ be the unit vector in this direction (unique up to a sign).

\vs

In both cases, by the discussion on the algebraic lemma and formula (2.4), we know that locally there will be orthonormal frame $\{ \xi_1, \ldots , \xi_4\}$ such that $A_{\xi_1}$ and $A_{\xi_2}$ are both of rank $2$ or less, and $A_{\xi_4}=JA_{\xi_3}$ has rank at least $6$. Furthermore, $E'=\mbox{span}\{ \xi_1, \xi_2\}$, as  the set of all normal directions in which the shape operator has rank $4$ or less, is uniquely determined. Also, if we restrict ourselves to a connected component $U$ in an open dense subset of of $M$, we may assume that in $U$ the orthonormal frame $\{ \xi_1, \xi_2\}$ of $E'$ is also uniquely determined, up to interchange and signs.

\vs

By letting $J\xi_3=\xi_4$ and $J\xi_4=-\xi_3$, we get an almost complex structure on $E$, the orthogonal complement of $E'$ in $N$. So to prove the main theorem, it suffices by Theorem 2 to find a developable ruling $L$ for $E$. This will follow from Codazzi equation (2.5) and a clever argument discovered by Dajczer and Gromoll in \cite{DG97}.

\vs

Consider $\eta=\xi_1$ or $\xi_2$. $A_{\eta}$ has rank $q\leq 2$. Denote by $\Delta_{\eta }$  the kernel of $A_{\eta}$ in $T$,  and by $\Delta_{\eta}^{\perp}$ its orthogonal complete in $T$. $\Delta_{\eta}^{\perp}$ is also the image space of $A_{\eta}$. First we claim the following:

\vvs

\noindent {\bf Claim:} {\em For either $\eta=\xi_1$ or $\eta=\xi_2$, the $E$-component of $ \nabla^{\perp }_v\eta$, denoted by $( \nabla^{\perp }_v\eta )^E$, is always $0$ for all $v\in \Delta_{\eta }$. That is, for any $v\in \Delta_{\eta }$, it holds}
\begin{eqnarray}
\langle \nabla^{\perp }_v\eta , \xi_3\rangle = \langle \nabla^{\perp }_v\eta , \xi_4\rangle  =0 .
\end{eqnarray}

\vs

To prove the claim, assume the contrary. Without loss of generality, we may assume that $\eta=\xi_1$ and there is a $v\in \Delta_{\eta }$ such that $\xi=( \nabla^{\perp }_v\eta)^E \neq 0$. By (2.5), since $A_{\eta }v=0$, we have
\begin{eqnarray}
A_{ \nabla^{\perp }_v\eta }u= A_{ \nabla^{\perp }_u\eta }v + \nabla_v(A_{\eta }u)+ A_{\eta }[u,v]
\end{eqnarray}
for any $u\in T$. Let $T_{\eta}=\{ u\in T\mid (\nabla^{\perp}_u\eta )^E=0\}$. Since $E$ is $2$-dimensional, the codimension of $T_{\eta}$ in $T$ is at most $2$.

\vs

Let $\{ e_1, \ldots , e_n\} $ be a frame of $V$ such that $\{ e_3, \ldots , e_n\}$ is a unitary frame of $V_0=\mbox{ker}(H) \cap \mbox{ker}(S')$ and is perpendicular to $\{ e_1, e_2\}$. We will also assume that $\{ e_{r+1}, \ldots , e_n\}$ is a unitary frame of $D\subseteq V$ corresponds to $\Delta_0$. So $\{ e_1, \ldots , e_r\}$ is a frame of $D^{\perp }$ corresponds to $\Delta_0^{\perp}\cong {\mathbb R}^{2r}$.

\vs

Let $W\subseteq T$ be the subspace corresponds to $V_0$ under the identification $V\cong T$. Note that $W\subseteq \Delta_{\xi_1}\cap \Delta_{\xi_2}$. Now consider the space $W'=W\cap \Delta_0^{\perp }$. Its real dimension is $2r-4\geq 6$ since $r\geq 5$, so the space $W''=W'\cap T_{\eta }$ is at least $4$ dimensional, as $T_{\eta }$ has codimension at most $2$ in $T$.

\vs

By (5.2), we know that for any $u\in W''$,  $A_{\xi}u$ is contained in the space
$$ \Delta_{\eta}^{\perp } +  \mbox{span} \{ A_{\xi_2}v\} ,$$
which has dimension at most $3$. So there will be $0\neq u_0\in W''$ such that $A_{\xi}u_0=0$. We have $A_{\xi_1}u_0=A_{\xi_2}u_0=0$ since $u_0\in W$. On the other hand, since $\xi\neq 0$, $\{ \xi , J\xi \}$  spans $E$, so by the fact that $A_{J\xi}=JA_{\xi}$, we get $A_{\eta '}u_0=0$ for any normal direction $\eta '$. This means that $\alpha_f(u_0, w)=0$ for any $w\in T$.

\vs

If we write $u_0=X+\overline{X}$ for (a unique) $X\in V$, then for any $Y\in V$, we have
$$ \alpha_f(u_0,Y) = S_{YX} + H_{Y\overline{X}} =0, \ \  \ \ \ \forall \ Y\in V.$$
Since $X\in W\subseteq \mbox{ker}(H)$, so we get $S_{YX}=0$ for any $Y$ thus $X\in \mbox{ker}(S)$ as well. This will force $X=0$ since we assumed that $u_0\in \Delta_0^{\perp }$. Thus $u_0=0$, a contradiction, and we have completed the proof of the claim.

\vvs

From the discussion in the algebraic lemma, we know that there will be local frame $\{ e_1, \ldots , e_n\}$ of $V$, such that $\{ e_3, \ldots , e_n\}$ is a unitary frame of $V_0$ and  is perpendicular to $\{ e_1, e_2\}$, and under this frame it holds
\begin{eqnarray*}
H^{\xi_1}&=& \mbox{diag}(1,0,0, \ldots , 0)\\
S^{\xi_1}&=& \mbox{diag}(a,0,0, \ldots , 0)\\
H^{\xi_2}&=& \mbox{diag}(0,\delta ,0, \ldots , 0)\\
S^{\xi_2}&=& \mbox{diag}(0,b,0, \ldots , 0)
\end{eqnarray*}
where $\delta =0$ or $1$, and $a$, $b$ are nonnegative. Write $e_i=\varepsilon_{2i-1}-\sqrt{-1}\varepsilon_{2i}$ for $1\leq i\leq n$, then under the real tangent frame $\{ \varepsilon_1, \ldots , \varepsilon_{2n}\}$, the first two shape forms are given by
\begin{eqnarray*}
A^{\xi_1}&=& \mbox{diag}(1\!+\!a,1\!-\!a,\ \ 0,\ \ \ 0; \ \ 0, \ldots , 0)\\
A^{\xi_2}&=& \mbox{diag}(\ \ 0,\ \ \ \   0,\  \  \delta \!+\!b, \delta \!-\!b; 0, \ldots , 0)
\end{eqnarray*}
Our goal is to show that there exists vector fields $v_1$ and $v_2$ on $M$ such that $L=\mbox{span}\{\xi_3-v_1, \xi_4-v_2\}$ satisfies $\langle \widetilde{\nabla }E', L\rangle =0$. That is, for any $i,j=1,2$,
$$ \langle \xi_{2+i} -v_i,\widetilde{\nabla }\xi_j\rangle =0 $$
or equivalently
\begin{eqnarray}
\langle \xi_{2+i} ,\nabla^{\perp }_u\xi_1 \rangle & = & \langle v_i ,A_{\xi_1}u \rangle \\
\langle \xi_{2+i} ,\nabla^{\perp }_u\xi_2 \rangle & = & \langle v_i ,A_{\xi_2}u \rangle
\end{eqnarray}
for each $i=1,2$ and any $u$ in $T$.

\vs

By the claim above, both sides of (5.3) are zero if $u$ is in the kernel space of $A_{\xi_1}$, which is spanned by $\varepsilon_3$ through $\varepsilon_{2n}$ and also $\varepsilon_2$ if $a=1$. So (5.3) just need to hold true for all $u\in \Delta_{\xi_1}^{\perp } = \mbox{Im}(A_{\xi_1 })$.

\vs

Similarly, both sides of (5.4) vanishes if $u$ is in the kernel of $A_{\xi_2}$, which is spanned by $\varepsilon_1$, $\varepsilon_2$, and $\varepsilon_5$ through $\varepsilon_{2n}$, and also $\varepsilon_4$ if $\delta =b$. So we just need (5.4) to hold true for all $u\in  \Delta_{\xi_2}^{\perp } = \mbox{Im}(A_{\xi_2 })$.

\vs

Since $\Delta_{\xi_1} + \Delta_{\xi_2}=T$, we must have $\Delta_{\xi_1}^{\perp} \cap \Delta_{\xi_2}^{\perp}=0$.  So we have direct sum decomposition
$$ T=(\Delta_{\xi_1} \cap \Delta_{\xi_2})\oplus \Delta_{\xi_1}^{\perp} \oplus \Delta_{\xi_2}^{\perp},$$
and $v_1$, $v_2$ can be uniquely determined in $\Delta_{\xi_1}^{\perp} \oplus \Delta_{\xi_2}^{\perp}$ by (5.3) and (5.4). But adding any element of $\Delta_{\xi_1} \cap \Delta_{\xi_2}$ onto $v_1$ or $v_2$ would not affect (5.3) or (5.4).  This establish the existence of developable ruling $L$ for $E$ and the proof the main theorem is complete in the non-minimal case.

\vvs

Next let us consider the minimal case, namely $H=0$ everywhere. By our previous discussion on the algebraic lemma, we know that either there exists a $2$-dimensional subspace $E'$ of $N$ in which the kernel of $S'$ has codimension at most $2$, and the orthogonal complement $E$ admits an almost complex structure $J$; or the entire normal bundle $N$ admits an almost complex structure $J$. In both cases, the almost complex structure is unique since no shape operator is allowed to vanish. We claim that $J$ is always admissible. This is automatic on any rank $2$ bundle, while in the case of $J$ on the rank four bundle $N$, we claim the following admissibility result:

\vvs

\noindent {\bf Theorem 4.} {\em Let $f:M^n \rightarrow {\mathbb R}^{2n+4}$ be a real K\"ahler submanifold such that there is an almost complex structure $J$ on $N$. Assume that no shape operator vanishes, and the rank $r\geq 2$ everywhere, then $J$ is admissible, namely, for any tangent vector $v$ and any normal field $\xi$, it holds}
\begin{eqnarray}
\nabla^{\perp }_vJ\xi =J\nabla^{\perp }_v\xi
\end{eqnarray}

\vs

Let us continue with our proof of the main theorem first, assuming that Theorem 4 is already established. In the case when $N$ itself is equipped with an almost complex structure $J$, Theorem 4 says that $J$ is admissible. So  by Theorem 1 in the previous section, we know that there is an isometric identification ${\mathbb R}^{2n+4} \cong {\mathbb C}^{n+2}$ under which $f$ becomes a holomorphic map. That is, $f: M^n\rightarrow {\mathbb C}^{n+2}$ is a holomorphic isometric embedding. Note that in this case,  any local piece of holomorphic hypersurface $Q^{n+1}$ containing (a piece of) $M^n$ would be a K\"ahler extension of $M$. So the conclusion of the main theorem holds in this case.

\vvs

\noindent {\em Proof of Theorem 4:} Let us choose a local orthonormal frame $\{ \xi_1, \ldots , \xi_4\}$ for the normal bundle $N$, so that $\xi_3=J\xi_1$ and $\xi_4=J\xi_2$. For any $1\leq \alpha , \beta \leq 4$, let us denote by $\phi_{\alpha \beta}$ the real $1$-form on $M$ given by $\langle \nabla^{\perp }\xi_{\alpha } , \xi_{\beta }\rangle $. Write the $4\times 4$ real, skew-symmetric matrix $\phi = (\phi_{\alpha \beta })$ in $2\times 2$ blocks:
$$ \phi = \left( \begin{array}{cc} \phi^1 & \phi^2 \\ - ^t\!\phi^2 & \phi^3 \end{array} \right) $$
It is easy to see that (5.5) is equivalent to $\phi^1=\phi^3$ and $\ ^t\phi^2=\phi^2$. Write
$$ (\phi^1-\phi^3) + \sqrt{-1}\ (\ ^t\!\phi^2-\phi^2) = \left( \begin{array}{cc} 0 & 1 \\ - 1 & 0 \end{array} \right) \lambda  ,$$
then it suffices to show that $\lambda =0$. Let $\{ e_1, \ldots , e_n\}$ be a unitary frame of $V$, and let  $\{ \varphi_1, \ldots , \varphi_n\}$ be its dual coframe of $(1,0)$-forms on $M$. Write $\langle \widetilde{\nabla } e_i, \xi_{\alpha }\rangle = \psi_i^{\alpha }$, then since $H=0$, each
$$ \psi_i^{\alpha }  = \sum_{j=1}^n S^{\alpha }_{ij}\varphi_j $$
is a $(1,0)$-form. Denote by $\psi^{\alpha }$ for the column vector $^t\!(\psi^{\alpha }_1, \ldots , \psi_n^{\alpha })$, and write
$$ \psi = (\psi ' ; \ \psi '') = (\psi^1, \psi^2; \ \psi^3,\psi^4).$$
By our choice of the normal frame, we have $\psi '' = -\sqrt{-1} \psi '$, therefore
 \begin{eqnarray}
 \psi =  (\psi ', -\sqrt{-1}\psi ').
 \end{eqnarray}
The connection matrix of $\widetilde{\nabla }$ under the frame $\{ e, \overline{e}, \xi \}$ is
$$ \tilde{\theta } = \left( \begin{array}{ccc} \theta & 0 & \psi \\ 0 & \overline{\theta } & \overline{\psi } \\ - ^t\!\overline{\psi } & - ^t\!\psi & \phi \end{array} \right) .$$
Applying (5.6) to the Codazzi equation $d\psi = \theta \psi + \psi \phi $, we get two equations. Multiplying the second equation by $\sqrt{-1}$, and take its difference with the first equation, we get
$$ \psi ' \left( \begin{array}{cc} 0 & 1 \\ - 1 & 0 \end{array} \right) \lambda  =0,$$
or equivalently, $\psi^1\wedge \lambda = \psi^2\wedge \lambda =0$. We claim that this will force $\lambda =0$, thus proving Theorem 3. Write $\lambda = \sum_{k} (a_k\varphi_k + b_k \overline{\varphi_k})$. The above equation on $\lambda $ means that for each $i$ and each $\alpha $,
$$ \sum_{j,k=1}^n   S^{\alpha }_{ij} a_k \varphi_j \wedge \varphi_k + \sum_{j,k=1}^n   S^{\alpha }_{ij}b_k \varphi_j \wedge \overline{\varphi_k} \ =\ 0.$$
The second part implies that
$ S_{ij}^{\alpha } b_k =0 $ for any $i,j,k$, thus $b_k=0$ for all $k$. The first part implies that
$S^{\alpha }_{ij}a_k=S^{\alpha }_{ik}a_j$ for any $\alpha $ and any $i,j,k$. Since $M$ has rank $r\geq 2$, there will be some combination $S=\sum t_{\alpha }S^{\alpha }$ so that $S$ is a complex symmetric matrix of rank at least $2$. Take a unitary matrix $P$ such that $^t\!P^{\!-\!1}SP^{\!-\!1}=D=\mbox{diag}(d_1, \ldots , d_n)$ is diagonal, with $d_1d_2\neq 0$. Then we have $S=\ ^t\!PDP$, and $S_{ij}a_k=S_{ik}a_j$ for any $i,j,k$ becomes
$$ d_lP_{lj}a_k= d_l P_{lk}a_j$$
for any $l,j,k$. Take $l=1$ and $2$, we notice that if $a_k$ are not all zero, then the first two rows of $P$ will be proportional, a contradiction. So we must have $a_k=0$ for all $k$. This completes the proof of Theorem 4. \qed

\vvs

Now we are left with the situation when there exists orthogonal decomposition $N=E'\oplus E$ such that $E$ is equipped with an almost complex structure $J$, and the kernel of $S'$ is at most $2$-dimensional. Here $S'$ is the $E'$-component of $S$. Write $V_0=\mbox{ker}(S')$ and denote by $k$ its codimension. $k$ is either $1$ or $2$. Let  $\{ \xi_1, \ldots , \xi_4\}$ be a local orthonormal frame of $N$ such that $\{ \xi_1, \xi_2\}$ is a frame of $E'$. We have $H=0$ and $S^{\xi_3} =  \sqrt{-1} S^{\xi_4}$.

\vs

By our previous discussion, we may exclude the possibility that $E'$ is also equipped with an almost complex structure. In other words, we may assume that
\begin{eqnarray}
S^{\xi_1} \neq \pm \sqrt{-1} S^{\xi_2}.
\end{eqnarray}
Also, the symmetry condition (2.3) holds for $S'$ as well. Our goal is to establish the existence of a developable ruling $L$ for $E$.

\vs

We will consider the case $k=2$ first. Let $\{ e_1, \ldots , e_n\}$ be a unitary frame of $V$, such that $\{ e_3, \ldots , e_n\}$ is a frame of $V_0=\mbox{ker}(S')$. As in the proof of Theorem 4, we will write $$\psi_i^{\alpha }=\langle \widetilde{\nabla}e_i, \xi_{\alpha }\rangle, \ \ \ \phi_{\alpha \beta } = \langle \nabla^{\perp } \xi_{\alpha }, \xi_{\beta }\rangle $$
and denote by $\theta$ the connection matrix of $M$ under $e$. We also let $\{ \varphi_1, \ldots , \varphi_n\}$ be the coframe of $(1,0)$-forms dual to $e$.

\vs

Note that since $\psi_i^{\alpha }=\sum_{j=1}^n S^{\alpha }_{ij}\varphi_j$, we have
$ \psi^3=\sqrt{-1}\psi^4$, where $\psi^{\alpha }$ stands for the $\alpha$-th column of $\psi$. Also,
$ \psi^1_i=\psi^2_i=0$ for each $i\geq 3$.

\vs

By the Codazzi equation $d\psi = \theta \psi + \psi \phi$, we get
\begin{eqnarray*}
d\psi^3 &=&\theta \psi^3 + \psi^1\phi_{13}+\psi^2\phi_{23}+\psi^4\phi_{43}\\
d\psi^4 &=&\theta \psi^4 + \psi^1\phi_{14}+\psi^2\phi_{24}+\psi^3\phi_{34}
\end{eqnarray*}
Multiplying $-\sqrt{-1}$ on the second line, and then add the result to the first line, we get from $ \psi^3=\sqrt{-1}\psi^4$ that
\begin{eqnarray} 0 \ = \ \psi^1(\phi_{13}-\sqrt{-1}\phi_{14})+\psi^2(\phi_{23}-\sqrt{-1}\phi_{24})
\end{eqnarray}
We will write $\sigma_1=\phi_{13}-\sqrt{-1}\phi_{14}$ and $\sigma_2=\phi_{23}-\sqrt{-1}\phi_{24}$. Write
\begin{eqnarray*}
 \psi^1_1 &=&  a\varphi_1 + b\varphi_2, \ \ \ \  \ \psi^1_2 \ = \ b\varphi_1 + c\varphi_2 \\
 \psi^2_1& =& a'\varphi_1 + b'\varphi_2, \ \ \ \ \psi^2_2 \ = \ b'\varphi_1 + c'\varphi_2
 \end{eqnarray*}
Since $S'$ also satisfies the symmetry condition (2.3), we have
\begin{eqnarray}
ac-b^2+a'c'-b'^2=0.
\end{eqnarray}
We first claim that both $\sigma_1$ and $\sigma_2$ must be linear combinations of $\varphi_1$ and $\varphi_2$. Assume otherwise, then by (5.8), we must have $\psi^1_1\wedge \psi^2_1=0$ and $\psi^1_2\wedge \psi^2_2=0$. So $(a,b)$ is proportional to $(a',b')$ and $(b,c)$ is proportional to $(b',c')$. The proportionality constants are also equal, so we have $S^1=\lambda S^2$ for some constant $\lambda$. Because $S'$ satisfies (2.3), we have $\lambda^2=-1$ since we assumed that $k=2$ here. So $S^1=\pm \sqrt{-1}S^2$, a contradiction to (5.7). So the claim must hold, and we can write
$$ \sigma_1 = \alpha \varphi_1 + \beta \varphi_2, \ \ \ \ \sigma_2 = \alpha ' \varphi_1 + \beta '\varphi_2.$$
The first two rows of (5.8) become
\begin{eqnarray}
a\beta -b\alpha + a'\beta ' - b' \alpha ' & = & 0 \\
b\beta - c\alpha + b'\beta ' - c' \alpha ' & = & 0
\end{eqnarray}
We claim that there exists $w_1$ and $w_2$ such that
\begin{eqnarray}
(\alpha , \beta ) & = & \ w_1(a,b)+ \ w_2(b,c) \\
(\alpha ', \beta ') & = & w_1(a',b')+ w_2(b',c')
\end{eqnarray}
hold simultaneously. First let us assume that $ac-b^2\neq 0$. Let $w_1$, $w_2$ be uniquely determined by (5.12), we have
\begin{eqnarray}
a\beta -b\alpha = w_2(ac-b^2), \ \ \ b\beta - c\alpha =w_1(b^2-ac).
\end{eqnarray}
If we write
$$ \delta_1=\alpha ' - (w_1 a'+w_2b'), \ \ \ \delta_2=\beta ' -(w_1b'+w_2c'),
$$
then we have
\begin{eqnarray*}
a'\beta ' - b' \alpha ' &=& w_2(a'c'-b'^2) +(a'\delta_2-b'\delta_1) \\
b'\beta ' -c'\alpha '&=&w_1(b'^2-a'c')+(b'\delta_2-c'\delta_1)
\end{eqnarray*}
Adding with (5.14), and using (5.9)-(5.11), we derive at
$$ \left( \begin{array}{cc} a' & b' \\ b' & c' \end{array} \right) \left[\begin{array}{c} \delta_2 \\ \!-\!\delta_1 \end{array} \right] = 0
$$
Since $a'c'-b'^2=-(ac-b^2)\neq 0$, we get $\delta_1=\delta_2=0$, so (5.12) and (5.13) hold.

\vs

If $ac-b^2=0$, then $a'c'-b'^2=0$ by (5.9). We claim that in this case $(a,b)$ cannot be proportional to $(a',b')$. Assume otherwise, say, $(a,b)=\lambda (a',b')$. Since $S^1$ and $S^2$ have zero determinants, we  have $(b,c)=\lambda (b',c')$ as well. So $S^1=\lambda S^2$, a contradiction to $k=2$, so the claim holds. Note that the claim means $\psi^1_1\wedge \psi^2_1\neq 0$. If we write $\psi^1_2=\lambda_1 \psi^1_1$ and $\psi^2_2 =\lambda_2 \psi^2_1$, then since $b=\lambda_1a$ and $b'=\lambda_2a'$, we know that $\lambda_1\neq \lambda_2$ by the above claim.

\vs

By (5.8), we have $\psi^1_1\sigma_1+\psi^2_1\sigma_2=0$ and $\lambda_1\psi^1_1\sigma_1+\lambda_2\psi^2_1\sigma_2=0$. Since $\psi^1_1\wedge \psi^2_1\neq 0$, the first equation implies that
$$ \sigma_1 = x\psi^1_1 + y\psi^2_1, \ \ \ \sigma_2 = y\psi^1_1 + z\psi^2_1 $$
for some scalar valued functions $x$, $y$, and $z$. Plug them into the second equation, we get $y(\lambda_1-\lambda_2)=0$, thus $y=0$. Take $w_2=(x-z)/(\lambda_2-\lambda_1)$ and $w_1=x-\lambda_1w_2$, we have $x=w_1+\lambda_1w_2$ and $z=w_1+\lambda_2w_2$, therefore
$$ \sigma_1=w_1\psi^1_1+w_2\psi^1_2, \ \ \ \sigma_2=w_1\psi^2_1+w_2\psi^2_2 $$
hold simultaneously. That is, (5.12) and (5.13) holds in this case as well.

\vs

Note that we have proved that, when $k=2$ and when $E'$ is not equipped with an almost complex structure,  there are scalar valued functions $w_1$ and $w_2$ such that $w=w_1e_1+w_2e_2$ satisfies $\sigma_1=\psi^1_w$ and $\sigma_2=\psi^2_w$, namely, for $\alpha =1$ and $2$, it holds that
$$ \langle \nabla^{\perp }\xi_{\alpha } , \ \xi_3\!-\!  \sqrt{-1} \xi_4 \rangle = \langle \widetilde{\nabla }w, \xi_{\alpha }\rangle .$$
If we write $w=-v_1+\sqrt{-1}v_2$, then the above just means that $\langle \widetilde{\nabla }E',L\rangle =0$ for the rank two subbundle $L$ in $T\oplus E$ spanned by $\{ \xi_3\!-\!v_1, \xi_4\!-\!v_2\}$. In other words, $L$ is a developable ruling of $E$. Thus by Theorem 2 we get a K\"ahler extension $h$ for $f$.  Note that since $w$ is a type $(1,0)$ vector, we have $v_2=Jv_1$ in this case. So $h$ is minimal by Theorem 3.

\vvs

Finally, let us consider the $k=1$ case, namely when $V_0=\mbox{ker}(S')$ has codimension one. Let $e=\{ e_1, \ldots , e_n\}$ be a unitary frame of $V$ so that $\{ e_2, \ldots , e_n\}$ is a frame of $V_0$. Let $\varphi$ be the dual coframe of $e$, and define $\psi$, $\phi$ as before. Then $\psi^3=\sqrt{-1}\psi^4$, and $\psi^1_i=\psi^2_i=0$ for all $i\geq 2$. Let us write $\psi^1_1=a\varphi_1$, $\psi^2_1=\lambda a\varphi_1$. Then $a\neq 0$, and $\lambda \neq \pm \sqrt{-1}$ since we have excluded the case where $S'$ admits an almost complex structure. By the Codazzi equation for $\psi^3$ and $\psi^4$, we again get
$$ \psi^1(\phi_{13}-\sqrt{-1}\phi_{14})+\psi^2(\phi_{23}-\sqrt{-1}\phi_{24})= \psi^1\sigma_1+\psi^2\sigma_2=0.$$
That is,
\begin{eqnarray}
\varphi_1(\sigma_1+\lambda \sigma_2)=0.
 \end{eqnarray}
 On the other hand,  since $\psi^4=-\sqrt{-1}\psi^3$, the Codazzi equation for $\psi^1$ and $\psi^2$ give
\begin{eqnarray*}
d\psi^1&=&\theta \psi^1 -\psi^2\phi_{12} - \psi^3\sigma_1 \\
d\psi^2&=&\theta \psi^2 +\psi^1\phi_{12} - \psi^3\sigma_2
\end{eqnarray*}
Now if we use the fact that $\psi^2=\lambda \psi^1$, we get $d\psi^2=d\lambda \wedge \psi^1+\lambda d\psi^1$, so the above two equations yield
$$ d\lambda \wedge \psi^1 = (1+\lambda^2)\psi^1\phi_{12} +\psi^3(\lambda \sigma_1-\sigma_2) $$
Looking at the $i$-th row of this equation, for any $i\geq 2$, we get
$$ \psi^3_i(\lambda \sigma_1-\sigma_2)=0, \ \ \ \forall \ 2\leq i\leq n.$$
If $\lambda \sigma_1-\sigma_2\neq 0$, then $\psi^3_i$ for all $2\leq i\leq n$ are multiples of $\lambda \sigma_1-\sigma_2$, which implies that the lower right $(n-1)\times (n-1)$ corner of $S^{\xi_3}$ will have rank at most $1$. This together with the fact that $S^{\xi_4}=-\sqrt{-1}S^{\xi_3}$ shows that $(S^{\xi_3}, S^{\xi_4})$, hence $S$, must have non-trivial kernel in $V_0$, since the dimension of $V_0$ is bigger than $2$. This contradicts the assumption that the rank of $M$ is at least $5$. So we must have
\begin{eqnarray}
\lambda \sigma_1-\sigma_2=0
\end{eqnarray}
Plug this into (5.15), and using the fact that $1+\lambda^2\neq 0$, we get $\varphi_1\sigma_1=0$, thus
$$ \sigma_1=w\psi^1_1, \ \ \  \sigma_2=\lambda \sigma_1 = w\psi^2_1 $$
for some $w$. Write $we_1=-v_1+\sqrt{-1}v_2$ for $v_1$ and $v_2$ real, we get
$$ \langle \widetilde{\nabla }E', \xi_3\!-\!v_1\rangle = \langle \widetilde{\nabla }E', \xi_4\!-\!v_2\rangle =0$$
That is, $L=\mbox{span}\{ \xi_3\!-\!v_1, \xi_4\!-\!v_2\}$ gives a developable ruling for $E$. Note that just like in the $k=2$ case, here we also have $v_2=Jv_1$, so $h$  is minimal by Theorem 3. This finishes the proof of the $k=1$ case, and the proof of the main theorem is now complete.

\vvs

Finally, let us remark that, in both the minimal and non-minimal cases, the K\"ahler extension is not necessarily unique, at least by the way we defined it, since one can add any vector fields in $\mbox{ker}(A_{E'})$ onto $v_1$, $v_2$, thus getting different developable rulings $L$. However, except in the case when $M^n$ is a complex submanifold of complex codimension $2$ in ${\mathbb C}^{n+2}$, there is always a `canonical' way to choose the developable ruling $L$, namely to take $L$ in such a way that $v_1$ and $v_2$ belong to
the orthogonal complement of $\mbox{ker}(A_{E'})$.  This uniqueness of canonical extensions might become important in the discussion of the global situations, namely, when $M$ is assumed to be complete.

\vvs

\vvs

\vvs

\vvs


\begin{thebibliography}{99}



\bibitem{BS1}   S. Brendle and R. Schoen, \emph{Manifolds with 1/4-pinched curvature are space forms},
  J. Amer. Math. Soc., \textbf{22} (2009), 287-307.

\bibitem{BS2}   S. Brendle and R. Schoen, \emph{Classification of manifolds with weakly 1/4-pinched
curvatures},  Acta Math., \textbf{200} (2008), 1-13.



\bibitem{Calabi}
E. Calabi, \emph{Isometric imbedding of complex manifolds}, Ann. Math., \textbf{58} (1953), no.1,1-23


\bibitem{D}
M. Dajczer, \emph{A characterization of complex hypersurfaces in ${\mathbb C}^m$}, Proc. Amer. Math. Soc.,
\textbf{105} (1989), no.2, 425-428.

\bibitem{D1}
M. Dajczer, \emph{Submanifolds and Isometric Immersions}, Publish or Perish, 1990.

\bibitem{DG}
M. Dajczer and D. Gromoll, \emph{Gauss parametrizations and rigidity aspects of submanifolds}, J.
Diff. Geom., \textbf{22}(1985), 1-12.

\bibitem{DG1}
M. Dajczer and D. Gromoll, \emph{Real K\"ahler submanifolds and uniqueness of the Gauss map},
J. Diff. Geom., \textbf{22}(1985), 13-28.

\bibitem{DG2}
M. Dajczer and D. Gromoll, \emph{The Weierstrass representation for complete minimal real  K\"ahler submanifolds }, Invent. Math., \textbf{119}(1985), 235-242.

\bibitem{DG97}
M. Dajczer and D. Gromoll, \emph{Real K\"ahler submanifolds in low codimension},
Diff. Geom. Appl., \textbf{7}(1997), 389-395.

\bibitem{DR}
M. Dajczer and L. Rodriquez, \emph{Rigidity of real K\"ahler submanifolds},
Duke Math. J., \textbf{53}(1986), 211-220.


\bibitem{FHZ}
L. Florit, W-S Hui, and F. Zheng, \emph{On real K\"ahler Euclidean submanifolds with non-negative Ricci curvature},   J. Eur. Math. Soc. (JEMS), \textbf{7} (2005), no.1, 1-11.


\bibitem{FZ1}     L. Florit and F. Zheng, \emph{On nonpositively curved Euclidean submanifolds: splitting
results}, Comm. Math. Helv., \textbf{74} (1999), no.1, 53-62.

\bibitem{FZ2}     L. Florit and F. Zheng, \emph{On nonpositively curved Euclidean submanifolds: splitting
results. II},  J. Reine Angew. Math., \textbf{508} (1999), 1-15.

\bibitem{FZ3}     L. Florit and F. Zheng, \emph{A local and global splitting result for real K\"ahler Euclidean
submanifolds}, Arch. Math.(Basel), \textbf{84} (2005), no.1, 88-95.

\bibitem{FZ-hypersurface} L. Florit and F. Zheng,  \emph{Complete real K\"ahler Euclidean hypersurfaces are cylinders},  Ann. Inst. Fourier (Grenoble), \textbf{57} (2007), no.1, 155-161.

\bibitem{FZ-codim2}  L. Florit and F. Zheng, \emph{Complete real K\"ahler submanifolds in codimension two},
Math. Z., \textbf{258} (2008),  no.2, 291-299.


\bibitem{MN}
F. Marques and A. Neves,  \emph{Min-Max theory and the Willmore conjecture}, preprint, arXiv:12026036

\bibitem{Nash}
J. Nash,  \emph{$C^1$ isometric embeddings}, Ann. Math. \textbf{60} (1954), no.3, 383-396.


\bibitem{Spivak}
V. Spivak,  \emph{A Comprehensive Introduction to Differential Geometry, Vol IV}, Publish or Perish, Berkeley, 1979.



\bibitem{Xu}    H.W. Xu, \emph{Recent developments in differentiable sphere theorem}, Fifth International
Congress of Chinese Mathematicians, AMS/IP, Studies in Advanced Math.,
Vol.51, 2012,  pp.415-430.


\bibitem{XG}    H.W. Xu and J.R. Gu, \emph{An optimal differentiable sphere theorem for complete
manifolds},  Math. Res. Lett., \textbf{17} (2010), 1111-1124.

\bibitem{XZ}    H.W. Xu and E.T. Zhao, \emph{Topological and differentiable sphere theorems for
complete submanifolds},    Comm. Anal. Geom., \textbf{17} (2009), 565-585.

\bibitem{Z1}     F. Zheng, \emph{Isometric embedding of K\"ahler manifolds with nonpositive sectional curvature}, Math. Ann., \textbf{304} (1996), no.4, 769-784.

\bibitem{Zheng}
F. Zheng, \emph{Complex Differential Geoemtry}, AMS and International Press, 2000.

\end{thebibliography}
\end{document}